\documentclass[a4paper]{amsart}

\usepackage{graphicx}

\newtheorem{theorem}{Theorem}[section]
\newtheorem{proposition}[theorem]{Proposition}
\newtheorem{corollary}[theorem]{Corollary}
\newtheorem{lemma}[theorem]{Lemma}

\newtheorem{Remark}[theorem]{Remark}
\newenvironment{remark}{\begin{Remark}\rm}{\end{Remark}}

\newcommand{\eps}{\varepsilon}
\newcommand{\la}{\lambda}
\let\cal=\mathcal
\let\Bbb=\mathbb

\begin{document}

\title[Tauberian Theorem for Laplace Transforms]{A Tauberian
Theorem for Laplace Transforms with  Pseudofunction Boundary
Behavior}

\dedicatory{To my young friend Larry Zalcman on his sixtieth
birthday}

\keywords{Distributions, Fourier transform, Laplace transform,
prime number theorem, pseudofunctions, Tauberian theory}

\subjclass[2000]{Primary: 40E05; Secondary:
11M45, 11N05, 42A38, 44A10, 46F20.}

\date{December 2, 2003}

\author{Jaap Korevaar}

\begin{abstract}
The prime number theorem provided the chief impulse for 
complex Tauberian theory, in which the boundary behavior of a
transform in the complex plane plays a crucial role. We
consider Laplace transforms of bounded functions. Our
Tauberian theorem does not allow first-order poles on the
imaginary axis, but any milder singularities, characterized by
pseudofunction boundary behavior, are permissible. In this
context we obtain a useful Tauberian theorem by exploiting
Newman's `contour method'.
\end{abstract}
\maketitle

\setcounter{section}{0}    
\section{Introduction} \label{sec:1}
In 1980 Don Newman \cite{N80} published a beautiful
proof for the prime number theorem (PNT) by complex analysis.
His vehicle was an old theorem of Ingham \cite{In35} involving
Dirichlet series, for which he found a clever proof by contour
integration. The method is easily adapted to give
Theorem \ref{the:1.1} for Laplace transforms; cf.\ the
author's paper \cite{Ko82} and Zagier \cite{Za97}. (Preprints
of these papers circulated shortly after Newman's article
appeared.) The contour method has recently been used in
numerous articles motivated by operator theory; see for
example Allan, O'Farrell and Ransford \cite{AFR87}, Arendt and
Batty \cite{AB88}, Batty \cite{Ba90}, and the book by Arendt,
Batty, Hieber and Neubrander \cite{ABHN01}. 

If one is interested only in a quick proof of the PNT, the
following result will suffice:
\begin{theorem} \label{the:1.1}
Let $a(\cdot)$ be (measurable and) bounded on
$[0,\infty)$, so that the Laplace transform
\begin{equation} \label{eq:1.1}
f(z)=\cal{L}a(z)=\int_0^\infty a(t)e^{-zt}dt,\quad z=x+iy,
\end{equation}
is well-defined and analytic throughout the open half-plane 
$\{x={\rm Re}\,z>0\}$. Suppose that $f(z)$ has an analytic
extension to the open interval $(-iB,iB)$ of the imaginary
axis. Then 
\begin{equation} \label{eq:1.2}
\limsup_{T\to\infty}\,\left|\int_0^T
a(t)dt-f(0)\right|\le\frac{2M}{B},
\quad\mbox{where}\;\;M=\sup_{t>0}\,|a(t)|.
\end{equation}
\end{theorem}
\begin{corollary} \label{cor:1.2}
If $a(\cdot)$ is bounded and $f=\cal{L}a$ extends analytically
to every point of the imaginary axis, the improper integral
\begin{equation} \label{eq:1.3}
\int_0^{\infty-}a(t)dt=\lim_{T\to\infty}\,\int_0^T
a(t)dt\;\;\mbox{exists and equals}\;\;f(0).
\end{equation}
\end{corollary}
Here the `Tauberian condition' that $a(\cdot)$ be bounded can
(in the real case) be replaced by boundedness from below.
However, this makes the proof more complicated; cf.\
\cite{Ko02} (section 9). In Section \ref{sec:2} we sketch how
to deduce the PNT.

Theorem \ref{the:1.1} and Corollary \ref{cor:1.2} are
contained in results of Karamata \cite{Ka34} (theorem B) and
Ingham \cite{In35} (theorem III), which were obtained by
Wiener's method \cite{Wi32}. They did not require that $f(z)$
can be extended analytically to every point of the imaginary
axis, but could get by with weaker boundary conditions. The
aim of the present paper is to reduce the boundary
requirements in Theorem \ref{the:1.1} to a minimum:
\begin{theorem} \label{the:1.3}
Let $a(\cdot)$ be bounded on $[0,\infty)$, so that 
the Laplace transform $f(z)=\cal{L}a(z)$, $z=x+iy$ is
analytic for $x={\rm Re}\,z>0$. Suppose that $f(x)$ tends
to a limit $f(0)$ as $x\searrow 0$ and that the quotient 
\begin{equation} \label{eq:1.4}
q(x+iy)=\frac{f(x+iy)-f(x)}{iy},\quad x>0,
\end{equation}
converges in distributional sense to a pseudofunction $q(iy)$
on the interval $\{-B<y<B\}$ as $x\searrow 0$. Then one has
inequality
$(\ref{eq:1.2})$.
\end{theorem}
Known sufficient conditions for (\ref{eq:1.2}) are 
uniform or $L^1$ convergence of $q(x+iy)$ to a limit
function $q(iy)$ on $(-B,B)$. The distributional conditions in
the Theorem require two things:

(i)$\;$ (convergence condition) that 
$$\int_{\Bbb R}q(x+iy)\phi(y)dy\;\;\;\mbox{should tend
to a limit}\;\;<q(iy),\phi(y)>$$ 
for every $C^\infty$ function $\phi$ with support in $(-B,B)$;

(ii) (pseudofunction condition) that $q(iy)$ be the restriction
to $(-B,B)$ of the distributional Fourier transform of a
function which tends to zero at $\pm\infty$. Cf. Sections
\ref{sec:4} and \ref{sec:5} below.

We remark that related distributional conditions received
inadequate treatment in \cite{Ko02} (Theorem 14.6). 
General background material on Tauberian theory can be found
in the forthcoming book \cite{Ko04}.

\setcounter{equation}{0} 
\section{From Corollary \ref{cor:1.2} to the Prime Number
Theorem}\label{sec:2}
Background material in number theory may be found in many
books; classics are Landau \cite{La09/53} and Hardy and Wright
\cite{HW79}. 

To obtain the PNT from Corollary \ref{cor:1.2} one may take
$a(t)$ equal to
\begin{equation} \label{eq:2.1}
b(t)=\frac{\psi(e^t)-[e^t]}{e^t}=e^{-t}\,
\sum_{1\le n\le e^t}(\Lambda(n)-1),
\end{equation}
where $\psi(v)=\sum_{n\le v}\Lambda(n)$ is Chebyshev's
function. The symbol $\Lambda(\cdot)$ stands for von
Mangoldt's function, which is given by the Dirichlet series
$$\sum_{n=1}^\infty\frac{\Lambda(n)}{n^w}
=-\frac{d}{dw}\log\zeta(w)=\frac{d}{dw}\sum_{p\,{\rm
prime}}\,\log(1-p^{-w})=\sum_{p\,{\rm
prime}}\,\frac{p^{-w}\log p}{1-p^{-w}}$$
when Re$\,w>1$.
It is elementary that $\psi(v)=\cal{O}(v)$, so that
$b(\cdot)$ is bounded. For Re$\,z>0$
\begin{eqnarray} \label{eq:2.2}
g(z) & = & 
\cal{L}b(z)=\int_0^\infty\{\psi(e^t)-[e^t]\}e^{-(z+1)t}dt
\nonumber \\ & = & \int_1^\infty\{\psi(v)-[v]\}v^{-z-2}dv
=\frac{1}{z+1}\int_{1-}^\infty v^{-z-1}d\{\psi(v)-[v]\}\\
& = & \frac{1}{z+1}\sum_1^\infty\frac{\Lambda(n)-1}{n^{z+1}}
=\frac{1}{z+1}\left(-\frac{\zeta'(z+1)}{\zeta(z+1)}
-\zeta(z+1)\right).\nonumber
\end{eqnarray}
The function $g(z)$ is analytic at every point of the line
$\{{\rm Re}\,z=0\}$. Indeed, $\zeta(w)$ is free of zeros on
the line $\{{\rm Re}\,w=1\}$ and the poles of
$-(\zeta'/\zeta)(w)$ and $\zeta(w)$ at the point $w=1$ cancel
each other. Conclusion: 
\begin{equation} \label{eq:2.3}
\int_0^{\infty-}b(t)dt
=\int_1^{\infty-}\frac{\psi(v)-[v]}{v^2}dv=g(0).
\end{equation}
By the monotonicity of $\psi$ this readily gives
$$\psi(v)\sim v\;\;\mbox{as}\;\;v\to\infty\quad\mbox{and}
\quad\sum_{n=1}^\infty\frac{\Lambda(n)-1}{n}=g(0).$$
The relation $\psi(v)\sim v$ is equivalent to the PNT:
$$\pi(v)\sim \frac{v}{\log v}\quad\mbox{as}\;\;v\to \infty.$$

\setcounter{equation}{0} \setcounter{figure}{2}
\section{An Auxiliary Result}
\label{sec:3}
We will prove Theorem \ref{the:1.1} but begin with a useful
preliminary form.
\begin{proposition} \label{prop:3.1}
Let $\sup_{t>0}\,|a(t)|=M<\infty$ and let the Laplace transform
\begin{equation} \label{eq:3.1}
f(z)=\cal{L}a(z),\quad z=x+iy,\;\;x>0,
\end{equation}
have an analytic extension to a neighborhood of the segment
$[-iR,iR]$ where $R>0$. Then for every number $T>0$,
\begin{eqnarray} & & \label{eq:3.2}
\Big|\int_0^T a(t)dt-f(0)\Big| \\
& & \le\frac{2M}{R}+\frac{|f(0)|}{eRT}+\frac{1}{2\pi}
\Big|\int_{-R}^R\{f(iy)-f(0)\}\left(\frac{1}{iy}+\frac{iy}{R^2}
\right)e^{iTy}dy\Big|.\nonumber
\end{eqnarray}
\end{proposition}
\begin{proof} Define
\begin{equation} \label{eq:3.3}
f_T(z)=\int_0^T a(t)e^{-zt}dt.
\end{equation}

(i) One begins with some simple estimates. For 
$x={\rm Re}\,z>0$,
\begin{equation} \label{eq:3.4}
|f_T(z)-f(z)|=\Big|\int_T^\infty a(t)e^{-zt}dt\Big|
\le M\int_T^\infty e^{-xt}dt=\frac{M}{x}e^{-Tx}.
\end{equation}
Similarly for $x={\rm Re}\,z<0$,
\begin{equation} \label{eq:3.5}
|f_T(z)|=\Big|\int_0^T a(t)e^{-zt}dt\Big|\le \int_0^T
Me^{-xt}dt <\frac{M}{|x|}e^{-Tx}.
\end{equation}

(ii) Let $\Gamma$ be the positively
oriented circle $C(0,R)=\{|z|=R\}$. We let $\Gamma_1$ be the
part of $\Gamma$ in the half-plane $\{x={\rm Re}\,z>0\}$,
$\Gamma_2$ the part in the half-plane $\{x<0\}$. Finally, let
$\sigma$ be the oriented segment of the imaginary axis from
$+iR$ to $-iR$ (Figure \ref{fig:3}). Observe that for
$z=x+iy\in\Gamma$, one has 
\begin{equation} \label{eq:3.6}
\frac{1}{z}+\frac{z}{R^2}=\frac{2x}{R^2}.
\end{equation}

By the hypotheses, the quotient $\{f(z)-f(0)\}/z$ is analytic
on the segment $\sigma$. Observe also that $f_T(z)$ is analytic
throughout the complex plane. Formulas
(\ref{eq:3.4})--(\ref{eq:3.6}) motivate the following
ingenious application of Cauchy's theorem and Cauchy's formula
due to Newman:
\begin{eqnarray} 
0 & = & \frac{1}{2\pi
i}\int_{\Gamma_1+\sigma}\frac{f(z)-f(0)}{z}dz
\nonumber \\ & = &  \label{eq:3.7} \frac{1}{2\pi
i}\int_{\Gamma_1+\sigma}\{f(z)-f(0)\}e^{Tz}
\left(\frac{1}{z}+\frac{z}{R^2}\right)dz, \\
f_T(0)-f(0) & = & \frac{1}{2\pi i}
\int_\Gamma\frac{f_T(z)-f(0)}{z}dz\nonumber\\ & = &
\frac{1}{2\pi i}\int_\Gamma\{f_T(z)-f(0)\}e^{Tz}\left(
\frac{1}{z}+\frac{z}{R^2}\right)dz. \label{eq:3.8}
\end{eqnarray}

\begin{figure}[htb] \label{fig:3}
$$\includegraphics[height=164bp]{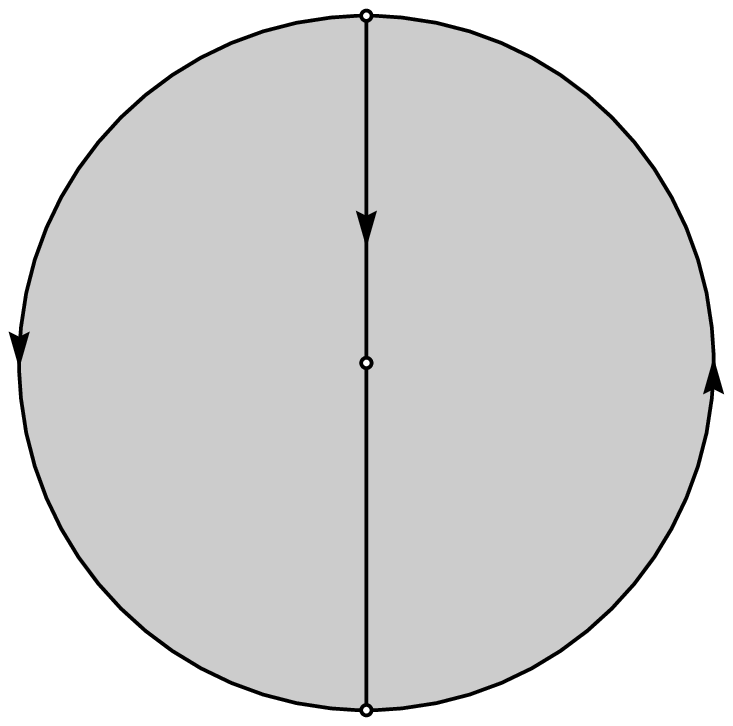}
\begin{picture}(0,164)(120,-72)
\put (20,10){$0$}
\put (105,10){$\Gamma_1$}
\put (-55,10){$\Gamma_2$}
\put (20,-70){$-iR$}
\put (20,85){$iR$}
\put (20,40){$\sigma$}
\end{picture}
$$
\caption{The paths of integration}
\end{figure}

\noindent Subtracting (\ref{eq:3.7}) from (\ref{eq:3.8}) and
rearranging the result, one  obtains the formula
\begin{eqnarray}  
f_T(0)-f(0) & = & 
\frac{1}{2\pi i}\int_{\Gamma_1}\{f_T(z)-f(z)\}e^{Tz}
\left(\frac{1}{z}+\frac{z}{R^2} \right)dz \nonumber\\
 & & +\, \frac{1}{2\pi i}\int_{\Gamma_2}\{f_T(z)-f(0)\}e^{Tz}
\left(\frac{1}{z}+\frac{z}{R^2}\right)dz \label{eq:3.9}\\ & &
-\, \frac{1}{2\pi i}\int_\sigma\{f(z)-f(0)\}e^{Tz}
\left(\frac{1}{z}+\frac{z}{R^2}\right)dz 
\nonumber \\ & = & I_1+I_2+I_3, \nonumber
\end{eqnarray}
say.

(iii) By (3.4) and (3.6) for $z\in\Gamma_1$, the integrand 
$f^*(z)$ in $I_1$ can be estimated as follows:
$$|f^*(z)| = \left|\{f_T(z)-f(z)\}e^{Tz}\left(\frac{1}{z}
+\frac{z}{R^2}\right)\right|\le \frac{M}{x}e^{-Tx}
e^{Tx}\frac{2x}{R^2}=\frac{2M}{R^2}.$$
Thus
\begin{equation} \label{eq:3.10}
|I_1|\le \frac{1}{2\pi}\int_{\Gamma_1}|f^*(z)||dz|
\le\frac{1}{2\pi}\frac{2M}{R^2}\pi R=\frac{M}{R}.
\end{equation}
For $z\in \Gamma_2$, where $|x|e^{Tx}\le 1/(eT)$, formulas
(3.5) and (3.6) imply the  estimate
\begin{equation} \label{eq:3.11}
|I_2|=\left|\frac{1}{2\pi
i}\int_{\Gamma_2}\{f_T(z)-f(0)\}e^{Tz}
\left(\frac{1}{z}+\frac{z}{R^2}\right)dz\right|
\le\frac{M}{R}+\frac{|f(0)|}{eRT}.
\end{equation}
Combination of (\ref{eq:3.3}) and
(\ref{eq:3.9})--(\ref{eq:3.11}) gives (\ref{eq:3.2}).
\end{proof}

\begin{proof}[Derivation of Theorem \ref{the:1.1}] 
Let $a$ and $f=\cal{L}a$ satisfy the hypotheses of Theorem
\ref{the:1.1}. Then we can apply
Proposition \ref{prop:3.1} for any $R\in (0,B)$. For the proof
of (\ref{eq:1.2}), one has to show that for any number
$\eps>0$, we can choose $T_0$ so large that the left-hand side
of (\ref{eq:3.2}) is bounded by $2(M/B)+\eps$ for all $T\ge
T_0$. To this end, choose $R$ so close to $B$ that
$2M/R < 2(M/B)+\eps/2$. In order to deal with the final term in
(\ref{eq:3.2}), or with 
\begin{equation} \label{eq:3.12}
I_3=\frac{1}{2\pi}\int_{-R}^R \{f(iy)-f(0)\}
\left(\frac{1}{iy}+\frac{iy}{R^2}\right)e^{iTy}dy,
\end{equation} 
one may apply integration by parts: $e^{iTy}dy=de^{iTy}/(iT)$, 
etc., or one may use the Riemann--Lebesgue lemma. 
Either method will show that for our $R$, 
\begin{equation} \label{eq:3.13}
I_3=I_3(R,T)\to 0\quad\mbox{as}\;\;T\to \infty.
\end{equation}
We now determine $T_0$ so large that 
$$\frac{|f(0)|}{eRT}+|I_3|<\eps/2,\quad\forall\,T\ge T_0.$$
Then by (\ref{eq:3.2})
$$\Big|\int_0^T a(t)dt-f(0)\Big|\le\frac{2M}{B}
+\eps,\quad\forall\,T\ge T_0.$$ 
\end{proof}

\setcounter{equation}{0}   
\section{Pseudofunction Boundary Behavior} \label{sec:4}
The preceding results may be refined with the aid of a
distributional approach. Motivated by operator theory,
Katznelson and Tzafriri \cite{KT86} used pseudofunctions to
strengthen the following theorem of Fatou \cite{Fa05},
\cite{Fa06}:
\begin{theorem} \label{the:4.1} Let the function
\begin{equation} \label{eq:4.1}
g(z)=\sum_{n=0}^\infty a_nz^n,\quad |z|<1,
\end{equation}
have an analytic continuation to (a neighborhood of) the point
$z=1$ on the unit circle $C(0,1)$. Suppose that the
coefficients satisfy the `Tauberian condition' $a_n\to 0$ as
$n\to\infty$.  Then the series $\sum_{n=0}^\infty a_n$
converges to $g(1)$.
\end{theorem}
The condition of analyticity at the point $z=1$ can be relaxed
in various ways. The most notable refinements in this
direction are due to M.\ Riesz and Ingham; cf.\ \cite{In35},
\cite{LaGa86}; another refinement is mentioned below.

The condition $a_n\to 0$ is the signature of {\it
pseudofunction} boundary behavior. In Fatou's theorem, and
for real $a_n$, it can be replaced by the one-sided condition
$\liminf a_n\ge 0$; cf.\ \cite{Ko54}, \cite{Ko02}. A
$2\pi$-periodic distribution $G(t)=\sum_{n\in{\Bbb
Z}}\,c_ne^{int}$ is called a pseudofunction if $c_n\to 0$ as
$n\to\pm\infty$. The latter condition first appeared
in Riemann's localization principle \cite{Ri92},
which Fatou used in the proof of his theorem. (A careful
discussion of the localization principle may be found in
\cite{Zy59}, item (5.7) in chapter 9.)

\medskip Let $g(z)$ as in (\ref{eq:4.1}) be any function
analytic in the unit disc. Among other things, Katznelson and
Tzafriri proved that pseudofunction boundary behavior of $g$ on
$C(0,1)\setminus\{z=1\}$, together with boundedness of the
sequence $\{s_n=\sum_{k=1}^n a_k\}$, implies that $a_n\to 0$. 
Their method can be used also for further relaxation of the
analyticity condition at the point $1$.
Knowing that $a_n\to 0$, it is enough for convergence of
$\sum a_n$ if $g$ in (\ref{eq:4.1}) is `weakly regular' at
the point $1$ in the following sense. For some constant which
may be called $g(1)$, the quotient
$$\frac{g(z)-g(1)}{z-1}$$ 
has pseudofunction boundary behavior at the point $z=1$ (more
precisely, in some angle $|\arg z|<\delta$); cf.\ 
\cite{Ko02}, \cite{Ko04}.

\medskip\noindent 
{\scshape Laplace Transforms} and related functions. Our aim is
to prove an extension of Theorem \ref{the:1.1} involving
pseudofunction boundary behavior of the Laplace
transform $f(z)=\cal{L}a(z)$. We begin with some general
remarks on tempered distributions, that is, continuous
linear functionals $F$ on the Schwartz space $\cal{S}$. The
`testing functions' $\phi\in\cal{S}$ include the $C^\infty$
functions with compact support. The result of applying $F$ to
$\phi$ is a bilinear functional, denoted by $<F,\phi>$.
Locally integrable functions $F_x(y)$ of at most polynomial
growth on $-\infty<y<\infty$ converge to a tempered
distribution $F(y)$ as $x\searrow 0$ if
$$\int_{\Bbb R}F_x(y)\phi(y)dy\to\,<F(y),\phi(y)>$$ 
for every function $\phi\in\cal{S}$. 

A tempered distribution $F$ on ${\Bbb R}$ is called a
{\it pseudomeasure} if it is the Fourier transform of a bounded
(measurable) function; it is called a {\it pseudofunction} if
it is the Fourier transform of a function which tends
to zero at $\pm\infty$. Reference: Katznelson \cite{Ka68/76}
(section 6.4). 

By the Riemann--Lebesgue theorem, every function in
$L^1({\Bbb R})$ is a pseudofunction. A nontrivial example of a
pseudomeasure on ${\Bbb R}$ is the distribution
$$\frac{1}{y-i0}=\lim_{x\searrow 0}\,\frac{i}{x+iy}=
\lim_{x\searrow 0}\,i\int_0^\infty e^{-xt}e^{-iyt}dt.$$ 
It is the Fourier transform of $i$ times the
Heaviside function, $1_+(t)$. Other examples are the Dirac
measure and the principal-value distribution, p.v.$\,(1/y)$.
In the case of boundary singularities, and roughly speaking,
first order poles correspond to pseudomeasures,
slightly milder singularities to pseudofunctions. 

Every pseudomeasure or pseudofunction $F$ on ${\Bbb R}$ can be
represented in the form 
\begin{equation} \label{eq:4.2}
F(y)=\lim_{x\searrow 0}\,\int_{\Bbb R}e^{-x|t|}b(t)e^{-iyt}dt,
\end{equation}
where $b(\cdot)$ is a bounded function, or a function
which tends to zero at $\pm\infty$, respectively. This
formula can be used to justify formal inversion of the order of
integration in some situations. An important consequence
is a {\it Riemann--Lebesgue lemma} for pseudofunctions $F$:
\begin{lemma} \label{lem:4.2}
For any pseudofunction $F$ on ${\Bbb R}$ and any testing
function $\phi$,
\begin{equation} \label{eq:4.3}
 <F(y),\phi(y)e^{iTy}>\;\to 0\quad\mbox{as}\;\;T\to\pm\infty.
\end{equation}
\end{lemma}
Indeed, by representation (\ref{eq:4.2}),
\begin{eqnarray*} 
<F(y),\phi(y)e^{iTy}> & = & \int_{\Bbb R}b(t)dt\int_{\Bbb
R}e^{-iyt}\phi(y)e^{iTy}dy\nonumber\\ & = & \int_{\Bbb
R}b(t)\hat\phi(t-T)dt\to 0\quad\mbox{as}\;\;T\to\pm\infty.
\end{eqnarray*}

\noindent{\scshape Products.} Let $F$ be a pseudomeasure
or pseudofunction as in (\ref{eq:4.2}) and let $\phi$ be a
testing function. Computing the Fourier transform of
$F(y)\phi(y)$, one finds that this product is the Fourier
transform of the convolution
$$\int_{\Bbb R}b(v-u)\hat\phi(u)/(2\pi)du. $$
For any other function $\Phi$ whose Fourier transform
$\hat\Phi(u)$ is $\cal{O}\{1/(u^2+1)\}$, the product $F\Phi$
may be {\it defined} as the Fourier transform of  
\begin{equation} \label{eq:4.4}
b^*(v)=\int_{\Bbb R}b(v-u)\hat\Phi(u)/(2\pi)du.
\end{equation}
With $F$, the product $F\Phi$ is again a pseudomeasure or
pseudofunction.

\setcounter{equation}{0}   
\section{Proof of Theorem \ref{the:1.3}} \label{sec:5}
Let $a(\cdot)$ and $f=\cal{L}a$ satisfy the hypotheses of the
Theorem. It is convenient to set $a(t)=0$ for $t<0$. Denoting
$\sup_{t>0}\,|a(t)|$ by $M$, taking $\eps>0$ and $0<R<B$, we
now apply Proposition \ref{prop:3.1} to $a(t)e^{-\eps t}$ and
$f(\eps+z)$ instead of $a(t)$ and $f(z)$. Thus we obtain the
inequality 
\begin{eqnarray} & & \label{eq:5.1}
\Big|\int_0^T a(t)e^{-\eps t}dt-f(\eps)\Big|
\\ & & \le\frac{2M}{R}+\frac{|f(\eps)|}{eRT}+\frac{1}{2\pi}
\Big|\int_{-R}^R\{f(\eps+iy)-f(\eps)\}\left(\frac{1}{iy}
+\frac{iy}{R^2}\right)e^{iTy}dy\Big|.\nonumber
\end{eqnarray}
To treat the final integral we set
\begin{equation} \label{eq:5.2}                      
\{f(\eps+iy)-f(\eps)\}\left(\frac{1}{iy}
+\frac{iy}{R^2}\right)=g_\eps(y).
\end{equation}
Let $\chi_R$ denote the characteristic function of the interval
$[-R,R]$. For any number $\la>0$ we let $\tau_\la$ denote a
`trapezoidal' testing function, that is, a $C^\infty$
function which is equal to $1$ on $[-\la,\la]$ and equal to
$0$ outside a suitable neighborhood of $[-\la,\la]$. The last
integral in (\ref{eq:5.1}) may then be written in
distributional notation as
\begin{equation} \label{eq:5.3}
I(T,\eps)=\,<g_\eps(y)\tau_R(y)\chi_R(y),e^{iTy}\tau_R(y)>.
\end{equation}
Here we take the support of $\tau_R$ inside $(-B,B)$. Then
by the hypotheses, $g_\eps(y)\tau_R(y)$ tends to the
pseudofunction 
$$g_0(y)\tau_R(y)=q(iy)(1-y^2/R^2)\tau_R(y)$$ 
as $\eps\searrow 0$; cf.\ (\ref{eq:1.4}). The question is
whether the integral $I(T,\eps)$ tends to the formal limit
$I(T,0)$. Multiplication by the cut-off function
$\chi_R(y)$ in (\ref{eq:5.3}) may cause problems! 

\medskip One may get around this difficulty by splitting the
integral $I(T,\eps)$. Choosing a trapezoidal function
$\tau_\mu$ with support in $(-R,R)$, we first consider the
relation
\begin{equation} \label{eq:5.4}
<g_\eps(y)\tau_\mu(y),e^{iTy}\tau_R(y)>\;\to\;
<g_0(y)\tau_\mu(y),e^{iTy}\tau_R(y)>\quad\mbox{as}\;\;
\eps\searrow 0.
\end{equation}  
By our Riemann--Lebesgue lemma \ref{lem:4.2}, the final
expression tends to zero as $T\to \infty$. 

Looking at (\ref{eq:5.3}), it remains to consider the `inner
product'
\begin{equation} \label{eq:5.5}
<g_\eps(y)\tau_R(y)\{1-\tau_\mu(y)\}\chi_R(y),
e^{iTy}\tau_R(y)>.
\end{equation} 
As $\eps\searrow 0$, the part of this expression which
comes from $f(\eps)$ tends to a trigonometric integral of an
integrable function,
$$\int_{-R}^R f(0)\left(\frac{1}{iy}
+\frac{iy}{R^2}\right)\{1-\tau_\mu(y)\}e^{iTy}dy.$$
The latter tends to zero as $T\to \infty$. From here on, we
focus on the constituent of the first factor in
(\ref{eq:5.5}) which involves $f(\eps+iy)$:
\begin{equation} \label{eq:5.6}
f(\eps+iy)\tau_R(y)\cdot\left(\frac{1}{iy}
+\frac{iy}{R^2}\right)\{1-\tau_\mu(y)\}\chi_R(y).
\end{equation}
The functions $f(\eps+iy)$ tend to the pseudomeasure 
$f(iy)=\hat a(y)$ as $\eps\searrow 0$, and by the hypothesis
about the quotient in (\ref{eq:1.4}), the restriction
of $f(iy)$ to $(-B,B)$ is equal to a pseudofunction. Hence the
product $f(iy)\tau_R(y)$, which by (\ref{eq:4.4})
is the Fourier transform of 
$$\frac{1}{2\pi}\int_{\Bbb R}a(v-u)\hat\tau_R(u)du,$$
{\it is} a pseudofunction. 

The functions $f(\eps+iy)\tau_R(y)$ are the Fourier
transforms of the functions $a(t)e^{-\eps t}$, which form a
uniformly bounded family. The factor
$$\Phi(y)=\left(\frac{1}{iy}+\frac{iy}{R^2}\right)
\{1-\tau_\mu(y)\}\chi_R(y),$$
which vanishes for $|y|\le\mu$ and for $|y|\ge R$, has 
Fourier transform $\hat\Phi(t)=\cal{O}\{1/(t^2+1)\}$. It
follows that the functions in (\ref{eq:5.6}) are
distributionally convergent. The limit $f(iy)\tau_R(y)\Phi(y)$
is a pseudofunction; cf.\ (\ref{eq:4.4}). The same will then be
true for the limit 
$$g_0(y)\tau_R(y)\{1-\tau_\mu(y)\}\chi_R(y)=\lim_{\eps\searrow
0}\,g_\eps(y)\tau_R(y)\{1-\tau_\mu(y)\}\chi_R(y)$$ 
of the functions in the first member of (\ref{eq:5.5}).
Combining the results, one concludes that the limit 
$I(T,0)$ of $I(T,\eps)$ can be written as an inner product
$$I(T,0)=\,<H(y),e^{iTy}\tau_R(y)>$$
involving a pseudofunction $H$, so that $I(T,0)\to 0$ as
$T\to\infty$.

To complete the proof of Theorem \ref{the:1.3} we return to
inequality (\ref{eq:5.1}). Letting $\eps$ go to zero one finds
that
\begin{equation} \label{eq:5.7}
\Big|\int_0^T a(t)dt-f(0)\Big|
\le\frac{2M}{R}+\frac{|f(0)|}{eRT}+\frac{1}{2\pi}
|I(T,0)|.
\end{equation}
Finally taking $T$ large and $R$ close to $B$, one obtains
the desired inequality (\ref{eq:1.2}). 
\begin{remark} \label{rem:5.1}
Related considerations show that one can introduce
pseudofunction boundary behavior in the statement of the
Wiener--Ikehara theorem \cite{Ik31}, \cite{Wi32}. One thus
obtains
\begin{theorem} \label{the:5.2}
Let $S(t)$ vanish for $t<0$, be nondecreasing, continuous from
the right and such that the Laplace--Stieltjes transform 
\begin{equation} \label{eq:5.8}
f(z)={\cal L}dS(z)=\int_{0-}^\infty e^{-zt}dS(t)
=z\int_0^\infty S(t)e^{-zt}dt,\quad z=x+iy,
\end{equation}
exists for $\mbox{\rm Re}\,z=x>1$.
Suppose that for some constant $A$, the analytic function
\begin{equation} \label{eq:5.9}
g(x+iy)=f(x+iy)-\frac{A}{x+iy-1},\quad x>1,       
\end{equation}
converges distributionally to a pseudofunction $g(1+iy)$ on
every finite interval $-B< y< B$ as $x\searrow 1$. Then
\begin{equation} \label{3eq:4.6}
e^{-t}S(t)\to A\quad \mbox{as}\;\; t\to \infty.
\end{equation}
\end{theorem}
\end{remark}

\bigskip

\noindent{\scshape Department of Mathematics, University of 
Amsterdam, \\
Plantage Muidergracht 24, 1018 TV Amsterdam, Netherlands}

\noindent{\it E-mail address}: {\tt korevaar@science.uva.nl}

\enddocument